\documentclass[11pt]{article}
\title{ODE, RDE and SDE Models of Cell Cycle Dynamics and Clustering in Yeast}

\author{
Erik M.~Boczko\\
Department of Biomedical Informatics\\
Vanderbilt University
\and
Tomas Gedeon\\
Department of Mathematics\\
Montana State University
\and
Chris C.~Stowers\\
Dow Chemical\\
Indianapolis, IN
\and
Todd R.~Young\footnote{Corresponding author, {\sf young@math.ohiou.edu}} \\
Department of Mathematics\\
Ohio University
}

\usepackage{amsmath}
\usepackage{epsfig}
\usepackage{graphics}
\usepackage{amsfonts}
\usepackage{amssymb}

\numberwithin{equation}{section}

\newtheorem{defn}{Definition}[section]

\newtheorem{thm}[defn]{Theorem}

\newtheorem{cor}[defn]{Corollary}

\newcommand{\D}[1]{{\mathbb#1}}% Doubled -Blackboard bold - caps only
\newcommand{\NN}{{\D{N}}}

\begin{document}

\maketitle

\begin{abstract}

Biologists have long observed periodic-like oxygen consumption oscillations
in yeast populations under certain conditions and several unsatisfactory
explanations for this phenomenon have been proposed. These ``autonomous oscillations''
have often appeared with periods that are nearly integer divisors
of the calculated doubling time of the culture.

We hypothesize that these oscillations could be caused by a weak form of
cell cycle synchronization that we call clustering.
We develop some novel ordinary differential equation models
of the cell cycle. For these models, and for random and stochastic perturbations,
we give both rigorous proofs and simulations showing that
both positive and negative growth rate feedback within the cell cycle
are possible agents that can cause clustering
of populations within the cell cycle. It occurs for a variety of models
and for a broad selection of parameter values. These results suggest that the
clustering phenomenon is robust and is likely to be observed in nature.
Since there are necessarily an integer number of clusters, clustering would
lead to periodic-like behavior with periods that are nearly integer divisors of
the period  of the cell cycle.

Related experiments have shown conclusively that cell cycle clustering
occurs in some oscillating yeast cultures.
\end{abstract}



\section{The Yeast Cell Cycle and Experimental Observations}

\begin{figure}[h!]
\centerline{\hbox{\psfig{figure=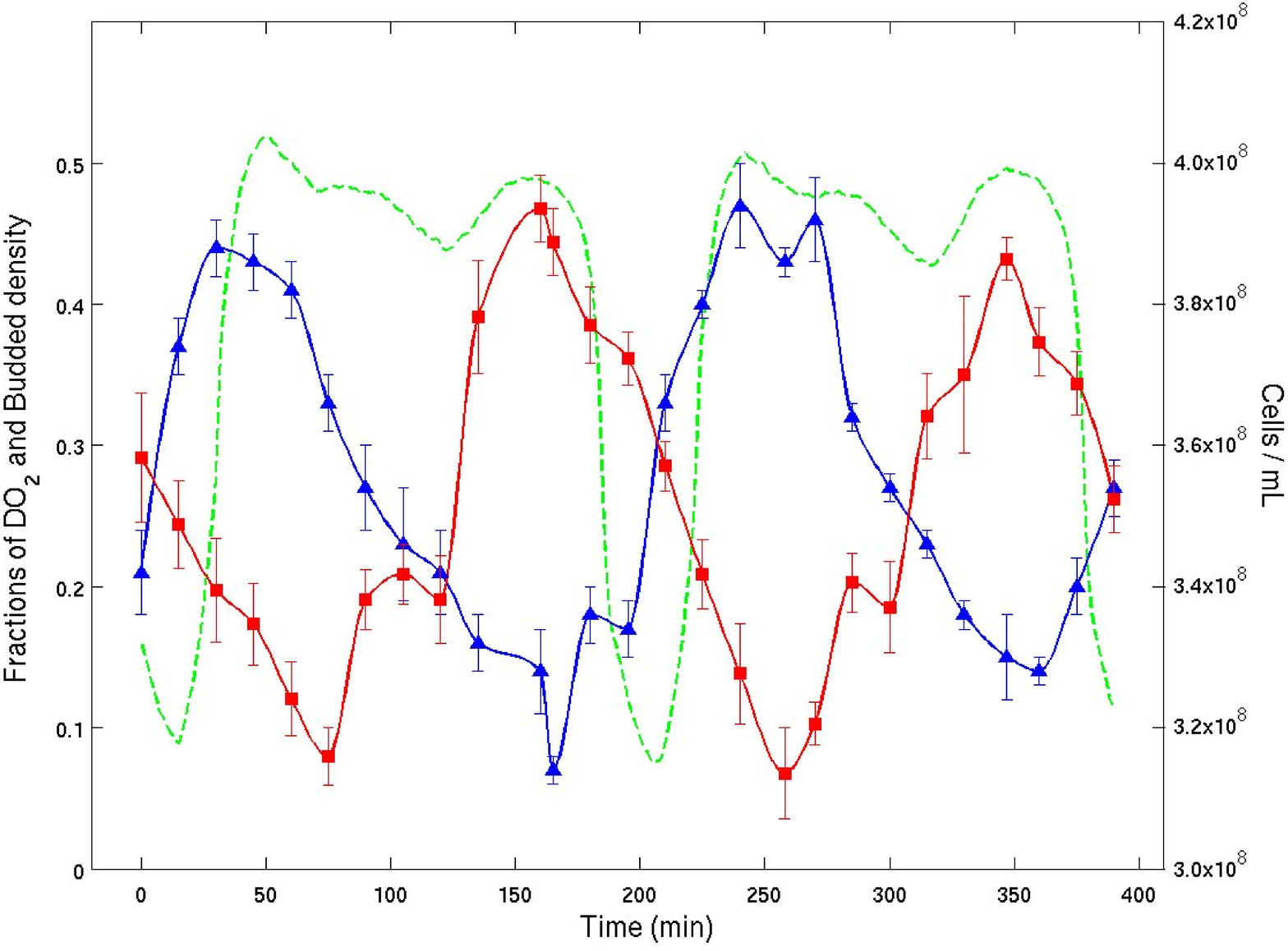,width=16cm,height=10cm}}}
\caption{Dissolved $O_2$ percentage (green), bud index percentage (blue)
and cell density (red). The plot shows clearly that the bud index (percentage
of cells with buds) and cell density  are both synchronized
with the oscillation in the level of dissolved $O_2$. Moreover, time periods of sharp
increases in bud index percentage show that a cluster of cell
(consisting of almost 50\% of the cells) is passing into S phase together. The doubling
time of this culture is 7.8 hours. Two oscillations per 7.8 hour period also is clearly
consistent with two clusters of cells.}
\label{ThePlot}
\end{figure}

Klevecz \cite{Kle1}, McKnight \cite{tu}, and others have observed
marked oscillations in the dissolved $O_2$ levels of various strains
of {\em Saccharomyces cerevisiae} (yeast) under certain conditions. For
instance in Figure~\ref{ThePlot} we show data from such an experiment
for the yeast strain Cen.PK~113 \cite{BBJRSY}.
It is seen in this experiment that dissolved $O_2$ levels oscillate between
about 5\% to 55\%, with a period of a little less than 4 hours. Since
the doubling time (as calculated from the dilution rate) for this culture was 7.8 hours
it is natural to suppose that the doubling time and $O_2$ oscillations are causally related.

Observations like these have been the subject of much speculation and
analysis \cite{boye,chenz,hense,henson,hjo,palk,tu,zhu}. Cultures exhibiting
$O_2$ oscillations, although planktonic,
are very dense; the average distance between cells is on the order of
1 cell diameter \cite{Kle1}. This suggests the possibility of cell to cell
or media signaling in the cultures. Klevecz and others suggested that
signaling is responsible for the $O_2$ oscillations by setting up
genomic synchronization. In the IFO0223 strain  they showed that
the compounds acetyaldehyde and $H_2S$ could be used to reset the
phase of the oscillations \cite{Kle1,Kle2}. Other
explanations of the oscillation involving various forms of signaling
have been proposed. For instance, \cite{hjo} suggested media signaling
as a mechanism for ``focusing'' populations.

During aerobic growth, the cell cycle of {\em Saccharomyces cerevisiae}
has four distinct phases, G1, S, G2 and M. The phase G1 begins
at cell division and is characterized by growth
of the cell. At the end of the G1 phase, which is thought to
be triggered by a volume milestone, the cell enters the S phase.
This is marked by the appearance of a bud and the beginning of replication
of the DNA. Replication continues throughout the phase. Once replication is
complete the cell enters a second growth phase G2. Most of the growth during
this period takes place in the bud. The M phase is marked by
narrowing  or ``necking'' of the connection between the original cell and the bud
and ends in cell division.

We propose here that dissolved $O_2$ oscillations in some strains could be due
to a phenomenon that we will call {\em clustering}. In general we
will define clustering loosely as significant groups of cells going
through cell cycle milestones at approximately the same time, i.e.~there
is a weak form of temporal synchronization. The clustering could cause
fluctuations in the $O_2$ dilution levels by clusters passing in and out of
high oxygen metabolism phases of the cell cycle, i.e. G1 phase.

In this article we study the implications of the hypothesis that the populations
of cells in the S-phase can  effect the growth rate of cells in the phase (denoted $R$
for responsive) before the S-phase, either positively or negatively.
We show rigorously and numerically in simple models that various forms of
signaling among cells, with either positive or negative feedback,
robustly leads to clustering of cells within the cell cycle.

We will prove mathematically
in idealized models that clustering occurs and that the
number of clusters formed depends on the widths $|S|$ and $|R|$
of the S-phase and the hypothesized region R.
Of key importance is the necessity that if clustering occurs, then
there are an integer number of clusters. This can provide the
basis for a periodic-looking behavior with a period which is approximately an
integer divisor of the cell-cycle length.
Models for which proofs are constructed include noisy or randomly
perturbed models. We use simulations to demonstrate
that the same type of behavior occurs in more detailed models.

In a separate experimental and analytical work reported elsewhere
\cite{BBJRSY} we show for the Cen.PK~113 strain in experiments similar to those
of \cite{tu} that bud index and cell density time-series
(see Figure~\ref{ThePlot}) prove definitively that clustering occurs.

In this paper we  focus primarily on models with the underlying assumption that
the cell cycles of different generations are identical, or nearly so.
Although this assumption is not justified in general yeast cultures, there
is some evidence for it in the cultures under study. In conditions where the
oscillations have been observed, the doubling time of the culture is
extended significantly, by a factor of about 4. It is known, and confirmed in  \cite{BBJRSY},
that when the cell cycle is extended, most of the extension occurs in the G1 phase.
It is also known that much of the differences between the cell cycles of
generations is due to differing lengths of the G1 phase (parent cells tend
to divide with larger volumes than their daughter cells and so need less time
to marshall sufficient resources to again undergo division). Thus by
extending the doubling time, the G1 phases of  all differing generations
are stretched comparably, which decreases the relative differences of
cell cycle lengths between generations. Age distributions
counts provide further evidence that the
assumption of equal cell cycle length for different
generations is approximately true in the cited experiments.

In addition to models with the assumptions of identical or nearly
identical cell cycles  we also
briefly treat models with  relative differences across generations
and cases where even larger differences could still lead to clustering.
In \cite{BBJRSY}, we also undertake a more careful study of a model
with stratified generations under various assumptions and using biologically
relevant values for all the parameters that can be determined. In all these
models clustering emerges as  a common and  robust phenomenon.

\section{Models of Cell Cycle Dynamics and Clustering}

\subsection{A general model for the cell cycle of an individual cell}

We will begin by defining a normalized logarithmic scale to represent the
cell cycle. Customarily, the cell cycle is delineated by volume milestones.
For a given cell, indexed by $i$, let $v_i(t)$ denote its volume.
It is usually assumed that the volume growth of a cell in a culture
is proportional to its volume, i.e.
$$
\frac{dv_i}{dt} = c_i(t) v_i
$$
where the relative growth rate $c_i(t)$ may depend on many factors,
such as available resources, chemical composition of the culture
substrate, etc. We will allow that the growth rate may also be influenced
by the state of the cell itself, thus the cell may react to environmental
factors in differing way at different stages of its cycle. The environmental
factors in turn may be influenced by the cells in the culture and their history.
Finally, the individual cells may have individual differences. Without any
loss of generality, we may include all of these factors in the rate $c_i(t)$.

Let  $V_{b,i}$ denote the volume of the cell at the
beginning of it cell cycle and $V_{d,i}$ its volume at division.
In order to have $x_i \in [0,1)$  we make the change of variables:
$$
     x_i = \frac{\ln(v_i / V_{b,i})}{\ln(V_{d,i} / V_{b,i})}.
$$
With this change birth has coordinate $x_i=0$ and division occurs at $x_i=1$.
We have then that $x_i$ satisfies
$$
\frac{dx_i}{dt} = d_i(t) \equiv  \frac{c_i(t)}{\ln(V_{d,i} / V_{b,i})}.
$$
Next we may rescale time so that
the average time span of a cell cycle in the culture is normalized to $1$, i.e.
scale by a factor $\bar{d}(t)$ which is the average (over all cells) growth
rate in a culture. In making this change of time we obtain
$$
\frac{dx_i}{dt} = \frac{d_i(t)}{\bar{d}(t)}.
$$
If the variation between individual cells is not too great then the right
hand side of this equation
is approximately $1$. Next, we distinguish between those influences
on the growth rate that are common to all the cells in the culture and those
due to individual differences. The common features we allow to depend on the
state, $x_i$ of the cell itself, as well as the conglomerate history of the
whole culture, which we denote by $\bar{x}$. Denote the common part of the growth
rate as $a$ and the individual part as $g_i$.
We can thus write the equation as:
\begin{equation}\label{rde}
\frac{dx_i}{dt} = a(t,x_i,\bar{x}) + g_i(t).
\end{equation}
The part of the growth rate due solely to individual differences in cells
is contained in the term $g_i$. In many applications this term will be
relatively small and ``random'' in the sense that differences are due
to details of the cell process that are far too minute and complex to model.
It contains both variations in the growth
rate and differences in $V_{b,i}$ and $V_{d,i}$. Thus we may
view (\ref{ode}) as a Random Differential Equation (RDE) \cite{arnold}.
A reasonable approximation of this equation in some circumstances is to
replace $g_i$ by a stochastic term:
\begin{equation}\label{sde}
 dx_i = a(t,x_i,\bar{x})\, dt + \sqrt{\sigma}\, dW_i,
\end{equation}
where $dW_i$ represents an independent noise term.
This Stochastic Differential Equation (SDE) must be interpreted in the
usual way as an Ito integral equation.
It is also reasonable under certain circumstances to consider (\ref{rde})
and (\ref{sde}) as perturbations of an Ordinary Differential Equation (ODE)
\begin{equation}\label{ode}
\frac{dx_i}{dt} = a(t,x_i,\bar{x}).
\end{equation}
This model allows for variation of the growth on state within the cell cycle
and the overall state of the system, but ignores differences between individual
cells. We will consider versions of this model in the rest of this section.
In the next two sections we will consider (\ref{rde}) and (\ref{sde}) as
perturbations of (\ref{ode}).

\subsection{Model of the culture}

First of all, we note that in a culture, cells are constantly
dividing, dying and perhaps being harvested. Thus the equations above must be applied
to a changing set of cells, indexed by a changing finite set of
positive integers $S(t) \subset \NN$. If the $i$-th cell dies or is harvested,
then $i$ is dropped from $S$. When a cell divides, one of
the new cells could continue being denoted by $i$, with $x_i$ resetting to
$0$ and the other (daughter cell)
would be identified by a new index. This  would be appropriate in budding yeast
where the mother and daughter cells are distinguishable.
The conglomerate history $\bar{x}$ would
contain each $x_i$ over its respective lifespan.

Now suppose that we are considering the unperturbed equation (\ref{ode}).
Since in this model there are no terms that differentiate the progression
of different cells, when a cell divides its two descendants will both start
at $x_i = 0$ and remain completely synchronized for the rest of their lifespans.
Thus tracking the evolution of both cells is redundant in this model and
we are inevitably led to consider a fixed set of cells. If the culture is
in a steady or periodic state, then there is also a probabilistic interpretation
of this simplification; given a
living cell at time $t$, the {\em expected number} of cells descended
from that cell at any later time $t+n$ is exactly $1$. Thus the original
model with a changing index set $S$ can be replaced by a fixed
system, each variable $x_i(t)$ representing the state of
its expected descendant at time $t$.

With a fixed index set, this model is easily amenable to numerical investigation
with currently available computer speeds. For instance, one can easily
investigate the behavior  of a conglomeration of 10,000 or more
cells over several cell cycles on a desktop workstation.
Further, under some additional assumptions on the dependence of $a$ and $g_i$,
we can investigate properties of the solutions of (\ref{ode}) and (\ref{rde})
rigorously as we show below.

\subsection{Modeling of the growth term $a$}

The standard biological assumption on $a$ is that
it does not depend on $x_i$ i.e. it is independent
of the cell's current state within the cycle. Rather
it is assumed to depend mostly on the nutrients available and
other environmental factors. With these assumptions $a = a(t,\bar{x})$.

As already described above, we hypothesize a more sophisticated form for $a$,
in which the location of cells in
cell cycle may influence the growth rate of themselves and other cells,
i.e.
$$
   a = a(t,x_i,\bar{x}).
$$
For example in budding yeast it has been hypothesized that cells in
the S phase may  influence the growth of cells in the
pre-budded G1 phase through some (unspecified) signaling mechanism.

In our models the population of
cells in the S-phase will be assumed to effect the growth of cells in
a preceding portion of the cell cycle which we denote by $R =[r,s]$, i.e.
a portion of the G1-phase. We will
consider a fixed finite population of $N$ cells, each of them
progressing via the following form of (\ref{ode})
\begin{equation}\label{ode-ide}
\frac{dx_i}{dt} = \begin{cases}
                 1 & \textrm{ if } \quad x_i \not\in R \\
                 1 +F(\#\{\textrm{cells in } S\}) & \textrm{ if } \quad x_i \in R.
                \end{cases}
\end{equation}
In this section we will consider two idealized (and discontinuous) forms of positive and
negative $F$ that we call the {\em advancing} and {\em blocking}.

By the {\em advancing model} we will mean that if the fraction of
cells in $S$ exceeds some threshold
$\tau$, then all cells in $R$ are instantaneously advanced
to the beginning of $S$, from which they resume normal growth.
This corresponds to a limit $F\to \infty $ in (\ref{ode-ide}) when
$\#\{\textrm{cells in } $S$\} \ge \tau$ and $F=0$ when $\#\{\textrm{cells in } $S$\} < \tau$.

By the {\em blocking model} we mean that if the fraction of
cells in $S$ equals or exceeds some threshold
$\tau$, then cells at $s$ are blocked from proceeding into
the S-phase, until the fraction of cells in $S$ drops below  the threshold.
This corresponds to choice $F=-1$ when $x_i = s$ and $\#\{\textrm{cells in } S \} \geq \tau$
and $F=0$ when $\#\{\textrm{cells in } $S$\} < \tau$. If more than $\tau$ cells have
accumulated at $s$, then all of those cells will enter $S$ together when
the fraction of cells drops below $\tau$. Note that periodic blocking at division was
considered by \cite{rotenberg}.

\subsection{Clustering}

Under either the advancing or blocking models, it is clear
that some cells may become synchronized. Consider the
advancing model; whenever the fraction of cells in $S$
exceeds the threshold $\tau$, then all cells in $R$
are instantaneously synchronized at $s$. There being no
mechanism in the model to subsequently differentiate
them, they will remain synchronized from then on.
Similarly, for the blockling model, during a period
when the threshold is exceeded, all cells arriving
at $s$ will be thereafter synchronized.
We will call a group of synchronized cells a {\em cluster}.
If the number of cells in a cluster is large enough to
exceed the threshold of the model, we call it a {\em critical cluster}.

Given an initial (discrete) distribution $\bar{x}(0)$ or
$w(s)= \sum_{i=1}^N \delta_{x_i(0)}(s)$ the evolution of such distribution under
either model will be called a {\em trajectory}  starting at $w(s)$ and denoted by
$\phi(t,w(s))$. It is clear that
different trajectories can merge.

A trajectory is called an {\em equilibrium} if it is stationary in the coordinate
frame moving with the speed $1$. In other words, $\phi(t,w(s))= w((x-t) \mod 1)$
for all $t \geq 0$.  A {\em periodic orbit } is a trajectory which is periodic in
the  moving frame. In such a case there exists $T$  with $\phi(T+t,w(x)) = \phi(t,w((x)))$
for all $t \in [0,T]$; if $T$ is the smallest number with this property, it is
called a {\em period.}

\subsection{Advancing model}

Denote by $\lfloor x \rfloor$ the largest integer  $\le x$ and
by $\lceil x \rceil$ the smallest integer  $ \ge x$.

\begin{thm}
Consider the advancing model.
\begin{description}
\item[A.]
If the initial  $w(x)$ exceeds the threshold on
all intervals of length $|S|$, except possibly inside the interval $R$, then the trajectory converges to a
periodic orbit.
\item[B.]
Any initial $w(x)$ that is below threshold on all intervals of length $|S|$, is an equilibrium point.
\item[C.] Every other initial condition converges to an equilibrium $e$ with a finite number of critical clusters, separated by voids of length at least $|R|+|S|$.
\end{description}
\end{thm}
\noindent
{\bf Proof:}
Let
\begin{equation}\label{qtx}
 q(t,x) = \int_{[x,x+|S|] \mod 1} \phi(t,w(z))dz.
\end{equation}
The function $q(t,x)$ is a form of local density of the discrete distribution of cells.
Observe that if $q(0,x) <\tau$ for all $x$ then this is a fixed point of the dynamics.
This corresponds to the case {\bf B} above.

The case {\bf A} corresponds to the case when
 $q(0,x) \geq \tau$ for all $x$, except for $x \in R$.
At $t=0$ all cells in $R$ are advanced  to the point $s$.
Since the $q(0,x) \geq \tau$ for all $x$ the cells arriving at $r$, the
start of the interval $R$, are immediately transferred to the point $s$. This
continues indefinitely. This trajectory is a periodic orbit in the moving frame.

Now we do the analysis of the case {\bf C}. Let $q(i,x), i = 0,1,\ldots$ be
the mass function after $i$ passes though the cell cycle. We decompose the domain
$I=[0,1]$ into intervals $I^i_1\cup J^i_1 \cup \ldots \cup J^i_{n-1} \cup I^i_k$ where
$q(i,x) \geq \tau $ for all $x \in I_j^i$ and $q(i,x) < \tau $ for all $x \in J_j^i$.
Set $I_j^i :=[ d_j^i, c_j^i ]$.
We compare the  intervals $I_j^i$ from iteration to iteration.
There are several possibilities that can happen to $I_j^i$ after a passage through a cell cycle:
\begin{enumerate}
\item The interval $I_j^i$ will shorten by $|R|$ and $I_j^{i+1} = [d_j^i +|R|, c_j^i].$
\item The intervals $I_j^i, \ldots, I_{j-k}^i$ may merge, if all intervals in between fall
within distance $|R|$ i.e. $c_j^i - c_{j-k}^i < |R|$;
\item If an interval  $I_{j-1}^i$ follows the  interval $I_j^i$ by a distance closer
then $|R| + |S|$ i.e.   $d_j^i -c_{j-1}^i \leq  |R| + |S|$, then the interval
 $I_{j-1}^i$, or a portion of that interval, will be promoted across $R$. The
resulting distance between the two new intervals will be less then $|S|$. After each subsequent
pass the distance between these two intervals will shorten by $|R|$ and in finite number of steps they will be within $|R|$ of each other and they will  merge. The number of steps this takes is uniformly bounded above by $l:=\frac{|S|}{|R|}+1$. This process  may result in a split of an interval $I_j^i$ into two intervals with the total mass conserved in the transaction. Notice that even though temporarily the number of intervals $I^{i+1}_j$ may increase over the number $I_j^i$,  after $k$ subsequent passes through $S$ the number of intervals $I^{i+k}_j$ is less or equal to number of intervals $I_j^{i}$.
\item The interval $I_j^i$ may shorten, if the preceding interval $I_{j+1}^i$
has critical mass in $S$ and therefore promotes part of the mass ahead of $c_j^i$ across $R$. Since a mass ahead of $c_j^i$ is lost, $c_{j}^{i+1} < c_j^i$.
The difference between this case and case (3.) is that here  $|R| + |S| \leq d_j^i -c_{j-1}^i \leq |R| + |S| + |S|$.
Observe, that the transferred cells will not form a new interval $I_j^{i+1}$,
since
they do not have a sufficient mass.
\end{enumerate}
In summary, the number of intervals $I_j^i$ is greater or equal to the number of intervals $I_j^{i+k}$ for every $i$
and the length of each intervals is a non-increasing function of $i$.
Finally, it is easy to see that the only interval $I_j^i$ that will not change
under the advance operator, is a singleton $I_j^i=[c_j^i,c_j^i]$. Since the population is finite, after a finite number of passes through $S$, each $I_j$ is a singleton. If two such singletons follow each other closer then $|R|+|S|$, by the step 3. above, they will merge in finite number of steps.
Eventually all remaining singletons $I_j$ are followed by a empty void of length $|R|+|S|$.

\hfill $\Box$

\begin{cor}
In the advancing model no more than $\lfloor(|R|+|S|)^{-1}\rfloor$
critical clusters can persist.
\end{cor}
\noindent
{\bf Proof:}
By the previous Theorem each critical cluster has to be followed by gap of length at least $|R|+|S|$.
\hfill $\Box$

%%%%%%%%%%%%%%%%%%%%%%%%%

\subsection{Blocking model}

\begin{thm}
In the blocking model no more than $\lceil |S|^{-1} \rceil$ critical clusters
can persist.
\end{thm}
\noindent
{\bf Proof:}
When two consecutive critical clusters pass through $S$, the second
cluster must wait at $s$ until the first cluster has passed $t$. Thus
the distance between them must then be at least $|S|$ and will remain
at least $|S|$ until the first cluster hits $r$. Then the first cluster
may have to wait to enter $S$, possibly decreasing the distance between the
two clusters in question.

Assume that the number $n$  of persistent clusters is constant along some trajectory.
Then either all these clusters are separated by a distance at least $|S|$, or there are
$n-1$ clusters separated by a distance $|S|$ and there is an additional cluster inside
$S$ whose distance to a waiting cluster at $s$ is smaller than $|S|$.

This minimal spacing implies that
$$
     (n-1) |S| \le 1.
$$
The result then follows.
\hfill $\Box$

We have shown that a cluster cannot persistently follow a
critical cluster closer than $|S|$ (except while the critical cluster is
being blocked at $s$). Also, it is clear that as many as
$\lceil |S|^{-1} \rceil$ critical clusters may persist if the
overall number $N$ of cell satisfies:
$$
N  > \tau  \lceil |S|^{-1} \rceil.
$$

Recall the definition of $q(t,x)$ in (\ref{qtx}).
\begin{thm}
In the blocking model, if the cells are initially distributed so that
the density $q(x,0)$ everywhere exceeds twice the threshold, then
exactly $\lceil |S|^{-1} \rceil$ clusters will develop and persist.
\end{thm}
\noindent
{\bf Proof:}
Cells at $s$ are initially stopped from entering $S$ until
the fraction of cells in $S$ drops below the threshold. Since
$ d(|S| - t) =\tau$, this
will occur at
$$
    t = |S| - \frac{\tau}{d}.
$$
At this time the cells that have clustered
at $s$ number $dt$ or $d|S| - \tau$. Since $d|S|$ is assumed
to be at least twice $\tau$, the threshold is again exceeded as the
cluster enters $S$. The threshold will continue to be exceeded until
this cluster leaves $S$. At this time a new cluster will have formed
at $s$ with volume $d|S|$ which is by assumption above $\tau$.
The second cluster is spaced exactly distance $|S|$ behind the first
cluster. At third cluster, etc. will form in the same way until the
first cluster returns to $s$. At that time there is either a cluster
in the interior of $S$ or at $t$. In the former case the first
cluster will wait less than time $|S|$ to enter $S$, ahead of when
the second cluster arrives. In the latter case it enters $S$ immediately.
This scenario will repeat itself indefinitely and the number of
critical clusters thus produced is an easy calculation.
\hfill $\Box$

\section{Small Perturbations of the Advancing and Blocking Models}

\subsection{Small Perturbations}

Next we consider {\em small perturbations} of the advancing
and blocking models by which we mean that the progression of each individual cell
is independently perturbed by a small term which is independent
of the other cells.
\begin{equation}\label{ode-perturbed}
\frac{dx_i}{dt} = \begin{cases}
                 1 + \epsilon \xi_i(x_i,t)& \textrm{ if } \quad x_i \not\in R \\
                 1 + F(\#\{\textrm{cells in } S\})
                       + \epsilon \xi_i(x_i,t)& \textrm{ if } \quad x_i \in R,
                \end{cases}
\end{equation}
where $|\xi_i(x_i,t)| \le 1$ and $\epsilon$ is small.
We will take as our definition of small that  $\epsilon$
is smaller than $\min(|S|,|R|)/6$. The effect of the perturbation
is to cause initially
synchronized cells to drift from each other, but
by no more than $2\epsilon$ or $\min(|S|,|R|)/3$ within a single cell cycle.
Note that under these assumptions (\ref{ode-perturbed}) can be considered as
a random differential equation with bounded noise. For general results about
RDE with bounded noise see \cite{homburg}.

Note that this model can viewed as a relaxation of the assumption that
cells are indistinguishable. It also can be considered as incorporating
random variations and noise.

With these perturbations, clusters will not remain precisely synchronized
as in the unperturbed models, but rather will spread out as the cells
proceed through the cell cycle. We now expand our definition of cluster
to include any group of at least $\tau$ cells that are within $|S|/3$
of each other in the cell cycle.

Below we will assume that the small perturbations act like noise in the
ways we define below. We say that perturbations satisfy the {\em maximum principle}
if
$$
   \max_{x} q(t_1,x) > \max_{x} q(t_2,x)
$$
and
$$
   \min_{x} q(t_1,x) < \min_{x} q(t_2,x)
$$
for any $t_1 < t_2$ that are not separated
by an activation or deactivation of blocking or advancing.
We will say that the perturbations are {\em diffusive} if:
\begin{itemize}
 \item the perturbations satisfy the maximum principle.
 \item synchronized cells will immediately de-synchronize.
\end{itemize}

\subsection{Perturbed Blocking Model}

\begin{thm}\label{perturbdelay}
In the blocking model with small diffusive perturbations such that $\epsilon < |S|/6$,
and an  initial distribution of cells that satisfies $q(0,x)> 2 \tau + 2 \epsilon$
for all $x$,
the system will form either $\lfloor |S|^{-1} \rfloor$, $\lceil |S|^{-1} \rceil$, or,
$\lceil |S|^{-1} \rceil + 1$ critical clusters within one cell-cycle period and
these clusters will persist indefinitely.
\end{thm}
\noindent
{\bf Proof:}
Denote $d = \min q(0,x)$.
Cells at $s$ are initially blocked from entering $S$ until
the fraction of cells in $S$ drops below the threshold. This
will occur at time $t_0$ no smaller than
$$
    t_0 = \frac{|S| - \frac{\tau}{d}}{1+\epsilon}.
$$
At this time the cells that have clustered
at $s$ number at least $t_0d$ or $(d|S| - \tau)/(1+\epsilon)$.
This group of cells we will call cluster 1.
By the assumptions, the threshold is again exceeded as the
cluster enters and crosses $S$. While crossing $S$ the cluster will de-synchronize,
but all the cells will remain in $S$ for at least $t_1 = |S|/(1+\epsilon)$
and no longer than $|S|/(1-\epsilon)$.
Thus while the first cluster crosses $S$ at least $t_1 d$ or $d|S|/(1+\epsilon)$
cell accumulate at $s$. This group we call the second cluster and from the
assumption $q(0,x) > 2 \tau + 2 \epsilon$ it is critical. Its crossing time $t_2$ will again be
bounded below by $|S|/(1+\epsilon)$ and above by $|S|/(1-\epsilon)$
and so this process of cluster formation will continue for as long as a density
of cells at least $d$ is arriving at $s$.

By a straight-forward calculation, the gap between any two adjacent clusters
formed in this process is bounded below by:
$$
     |S| \frac{1-\epsilon}{1+\epsilon}
$$
and obviously the gap is bounded above by $|S|$. If the first of two adjacent clusters
travel at the minimal speed $1-\epsilon$ while outside of $S$, all the cells in the
first cluster must arrive again at $s$ no later than:
$$
    t^* = \frac{1}{1-\epsilon} - \frac{|S|}{1+\epsilon},
$$
after cluster 2 leaves $S$.
During the same time period the second of the adjacent clusters can travel at a rate at
most $1+\epsilon$ and thus can travel a distance of no more than
\begin{equation*}
  d^* = t^* (1+\epsilon) = \frac{1+\epsilon}{1-\epsilon} - |S| < 1.
\end{equation*}
Therefore the entirety of the first cluster must reach $s$ before any of the second cluster does so. The gap at that time is at least
$$
   1-d^* = |S| - \frac{2 \epsilon}{1+\epsilon} > \frac{2 |S|}{3}.
$$

Next consider the relative motion between the cells that began in $S$ and the cells in cluster 1.
Note that the initial gap will be at least $|S|/2$. Since $\epsilon < |S|/4$, the last of the
cells will reach $s$ before the first cluster returns.

When cluster 2 leaves $S$, it is easily shown that the gap between cluster 1 and 2 is
no wider than
$$
     |S| \frac{1+\epsilon}{1-\epsilon}
                < |S|(1+2\epsilon).
$$
At the same time the gap between cluster 1 and 3 is no more than the gap between cluster
1 and 2 plus $|S|$. That is,
$$
      2 |S| (1+\epsilon).
$$
When cluster 3 reaches $t$, the gap between clusters 1 and 3 is no more than
$$
     2 |S| (1+\epsilon) + |S| \frac{1+\epsilon}{1-\epsilon} - |S| < 2(1 + 2 \epsilon)|S|.
$$
By induction one can show that the first cluster can be no more that
$$
    (j-1)(1+2\epsilon)|S|
$$
ahead of the $j$-th cluster when the $j$ cluster is at $t$. From this we can show
that at least $\lfloor |S|^{-1} \rfloor$ clusters will form.

On the other hand, when cluster 2 leaves $|S|$ the gap between clusters 1 and 2 must
be at least
$$
  \frac{2}{1+\epsilon}|S|
$$
and cluster 1 must be at least
$$
     \frac{j-1-\epsilon}{1+\epsilon} |S|
$$
ahead of the $j$-cluster. It is clear then that $j$ must be less than
$\lceil |S|^{-1} \rceil + 1$.

Now to show that the clusters persist.
Clearly, while cluster 1 is passing through the cycle the first time,
there is always a critical cluster in $S$, at least until the cells that
began in $S$ again reach $S$. Note that the cells which were inside $S$
when the first cluster was released from $s$ are critical in number. We will
call this group of cells the tail.
When the tail begins to arrive back at $S$, its leading edge can be no more than
$$
     \frac{|S|(1-|S|)}{1-\epsilon}
$$
ahead of cluster 1. At that moment there will either be (1) a critical in the interior
of $S$, or, (2) there are critical clusters at both ends of $S$.
Let us consider case (2) first. In this case the cluster at the beginning of $S$ will
be critical and as it enters $S$ the cells from tail will be blocked at $s$. Since the
original width of the tail is less than $|S|/2$, its width when it reaches $S$ will
be less than $2|S|/3$ and it will all reach $S$ before the cluster ahead of it leaves
$S$. Now two things can happen, both of which continue the process, either all or part
of cluster 1 reaches $S$ while the tail is blocked there, or it will be blocked by the
tail. In the former situation part or all of cluster 1 will merge with the tail and
be blocked and in the later the entirely of cluster will be blocked by the tail. Either way
cluster 1 is blocked and will be inside $S$ when cluster 2 arrives back. Now in case (1)
all or part of the tail will reach $S$ and be blocked while the previous cluster clears
$S$. It will merge with other cells blocked at $S$ and they will be a critical cluster
when they are released. Thus when cluster 1 arrives at $S$, the last cluster including
part or all of the tail will still be in $S$. Again in this case cluster 1 is blocked
before reentering $S$ and will be inside $S$ when cluster 2 arrives. Now in either case
cluster 1 is not only blocked, but resynchronized as it is blocked. The same happens for
the subsequent clusters and the process continues.
\hfill $\Box$

\subsection{Perturbed Advancing Model}
\label{pertadv}

First we note the lack of rigorous results for this model. We are not
able to prove anything similar to Theorem~\ref{perturbdelay}. This
lack of results, however is suggestive that the advancing model does
not form clusters as naturally as does the blocking model.

For this model we can make the following observation:
{\em In the advancing model with small random perturbations and $|R|=|S|$,
 an initial single critical cluster will have ``forward leakage.''}
Wherever a cluster begins in the cell cycle, when its edge reaches $s$
it will have positive width. As some of the  cells enter $S$ the threshold
will be exceeded, and the rest of the cells will be promoted to
$s$. The cluster will then proceed around the cell cycle. When its
edge reaches $S$ again, its width will be larger and the leading edge
flatter than during the previous cycle. Thus the front edge will reach deeper
into $S$ before the advance takes place. Thus the front edge of the
cluster will grow wider
and wider on each pass until some of the cells are able to pass through
$S$ completely before advance occurs. At this point these cells have essentially
escaped from the cluster.

Note that the cells which leak forward may eventually slow down and drop
back into the cluster, or, if they continue at a faster rate than the cluster will
eventually be caught in the next cluster (or into the original cluster if it is
the only one.) This implies that there may be steady states which have persistent
clustering, but the cells in the clusters may not actually be synchronized in the
strict sense.

\section{Linear feedback models and simulations}

\subsection{Graduated but unstratified model.}

Next we consider more realistic models. Rather than strict advancing
and blocking we consider the possibility that the number of cells in
$S$ either act to slow down or speed up cell growth in the preceding region $R$.
The population of cells in the
S-phase is hypothesized to effect the growth of cells in
a preceding phase $R =[r,s]$ via a continuous function. This influence
will be assumed to be a function of the number of cells
in the $S$ phase. Thus we have the following equations of
motion:
\begin{equation}\label{sde_SR}
dx_i = \begin{cases}
             dt + \sqrt{\sigma} \, dW_i & \textrm{ if } \quad x_i \not\in $R$ \\
             ( 1 + F(\#\{\textrm{cells in } $S$\}) ) \, dt
                          + \sqrt{\sigma} \, dW_i & \textrm{ if } \quad x_i \in $R$.
                \end{cases}
\end{equation}
Note that this form is quite general since we might consider any
functional dependence $F$ and $R$ could be large or small.

We report simulations of this model where we will assume $F$ to be linear
function  $F(x) = ax$, with $a<0$ for the negative feedback model
and $a>0$ for the positive model.
For example a nonlinear model would be obtained if $F$ was a sigmoidal
Hill function $F(x) = \frac{A}{1+ (x/\tau)^n}$. In the negative case in the limit of Hill coefficient $n \to \infty$ we obtain a strict blocking model described  in section 2.

In the simulations we also incorporate a small diffusive term (white noise)
with variation $\sigma$. The program tracked the trajectories of 5000 individual
cells, initially uniformly distributed, through 20 unperturbed cell cycles
(20 units of time in equation \ref{sde_SR}). The simulation used an Euler method
to integrate the SDE. In each of Figures 1-4 we show
(a) a histogram of the final distribution of cells within the cell cycle
and (b) a time-series of the fraction of cells inside $S$ over the final
two periods.

\begin{figure}[h!]
\centerline{\hbox{{\psfig{figure=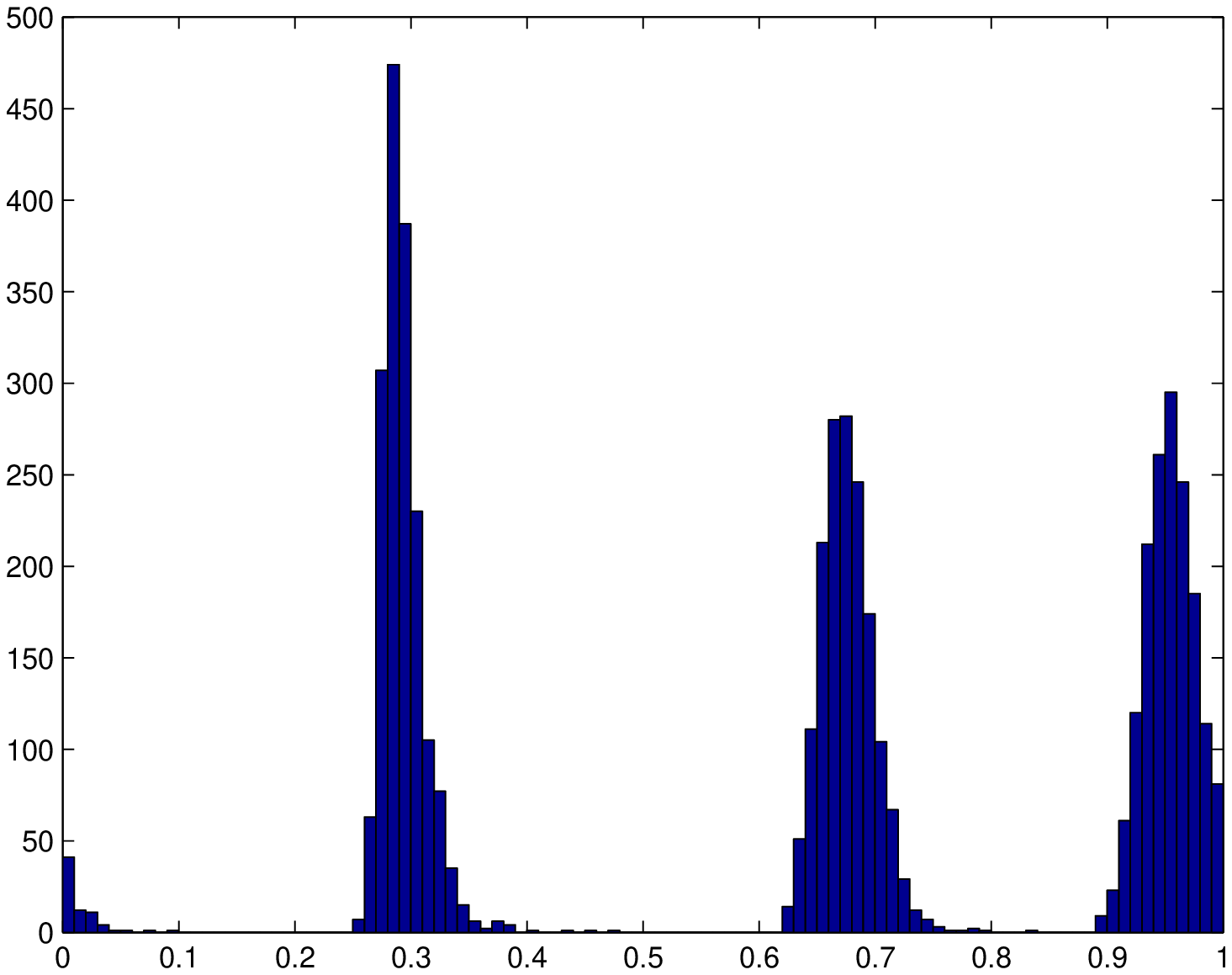,height=5cm,width=5cm}}
\hspace*{.1in}{\psfig{figure=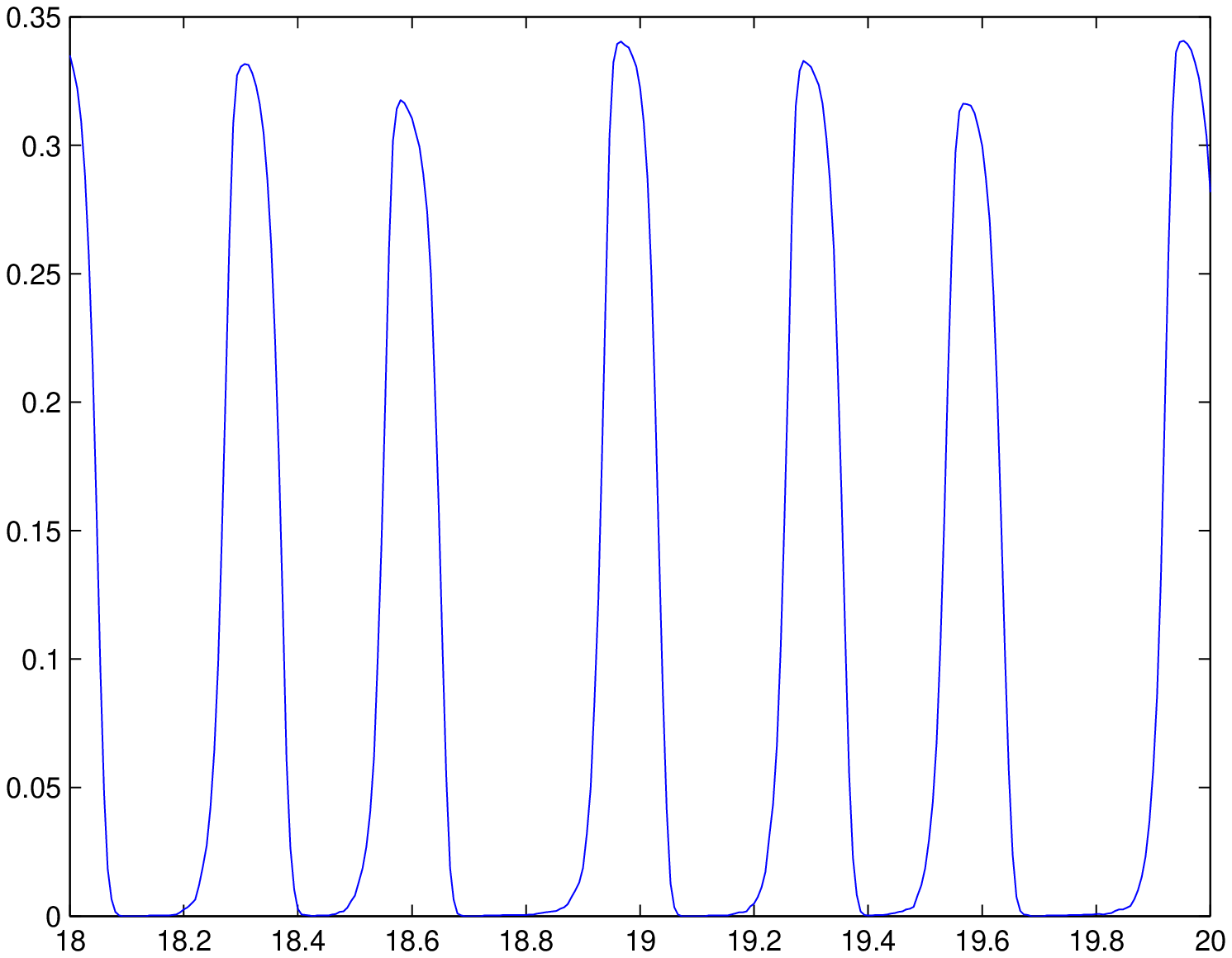,height=5cm,width=5cm}}}}
\caption{Simulations of the graduated model with positive feedback. Here
$R=[.1,.2]$, $S=[.2,.3]$, $\sigma = .02$ and the
linear factor was 5. (a) Histogram of the final distribution of
cells within the cell cycles. (b) Time-series of the final two time frames. One unit of time corresponds to one unperturbed cell cycle.}
\label{figpitch1}
\end{figure}

\begin{figure}[h!]
\centerline{\hbox{{\psfig{figure=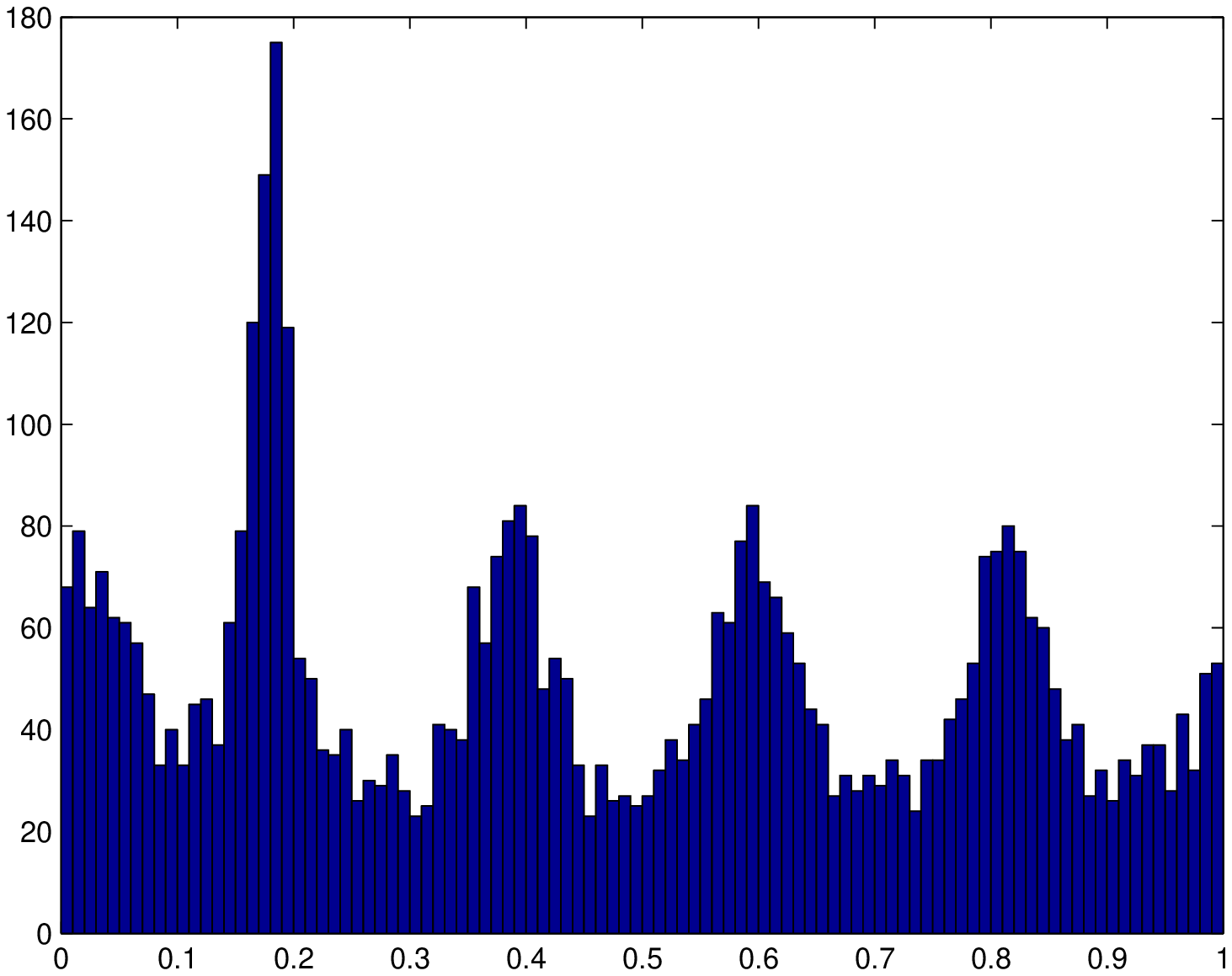,height=5cm,width=5cm}}
\hspace*{.1in}{\psfig{figure=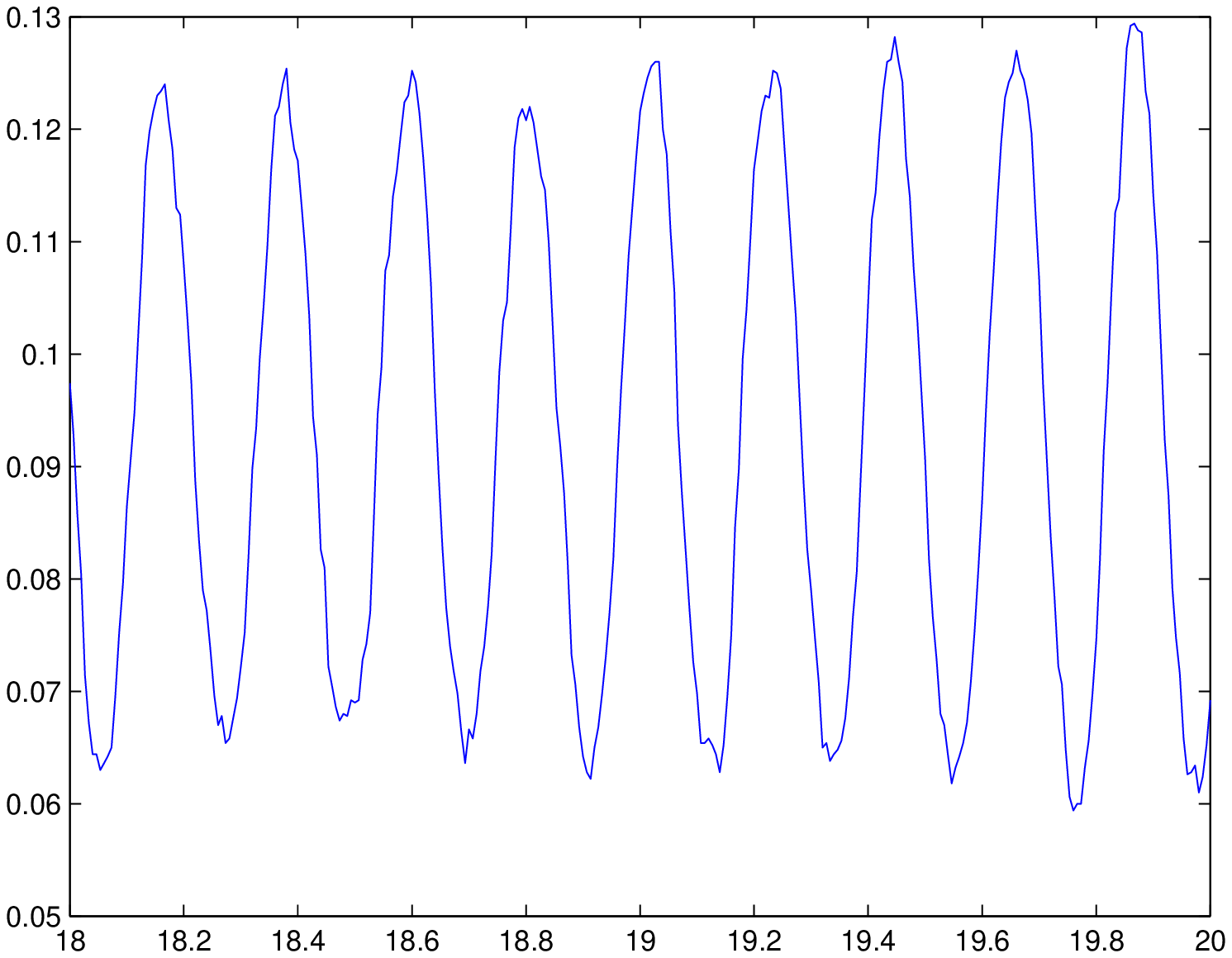,height=5cm,width=5cm}}}}
\caption{Simulations of the graduated model with negative feedback. Here
$R=[.1,.2]$, $S=[.2,.3]$, $\sigma = .02$ and the
linear inhibition factor was $-5$. (a) Histogram of the final distribution of
cells. (b) Time-series of the final two time frames. }
\label{figpitch2}
\end{figure}

\begin{figure}[h!]
\centerline{\hbox{{\psfig{figure=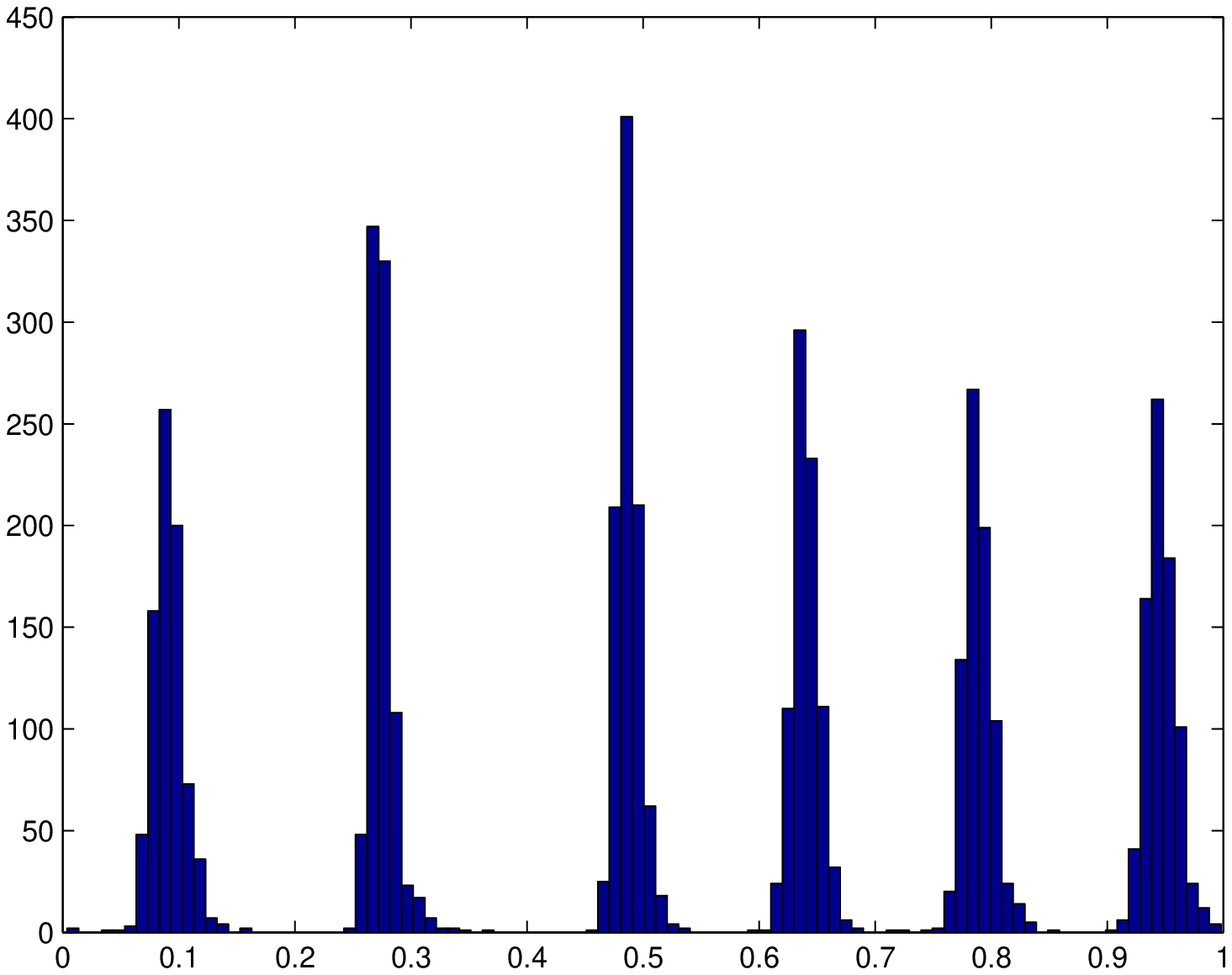,height=5cm,width=5cm}}
\hspace*{.1in}{\psfig{figure=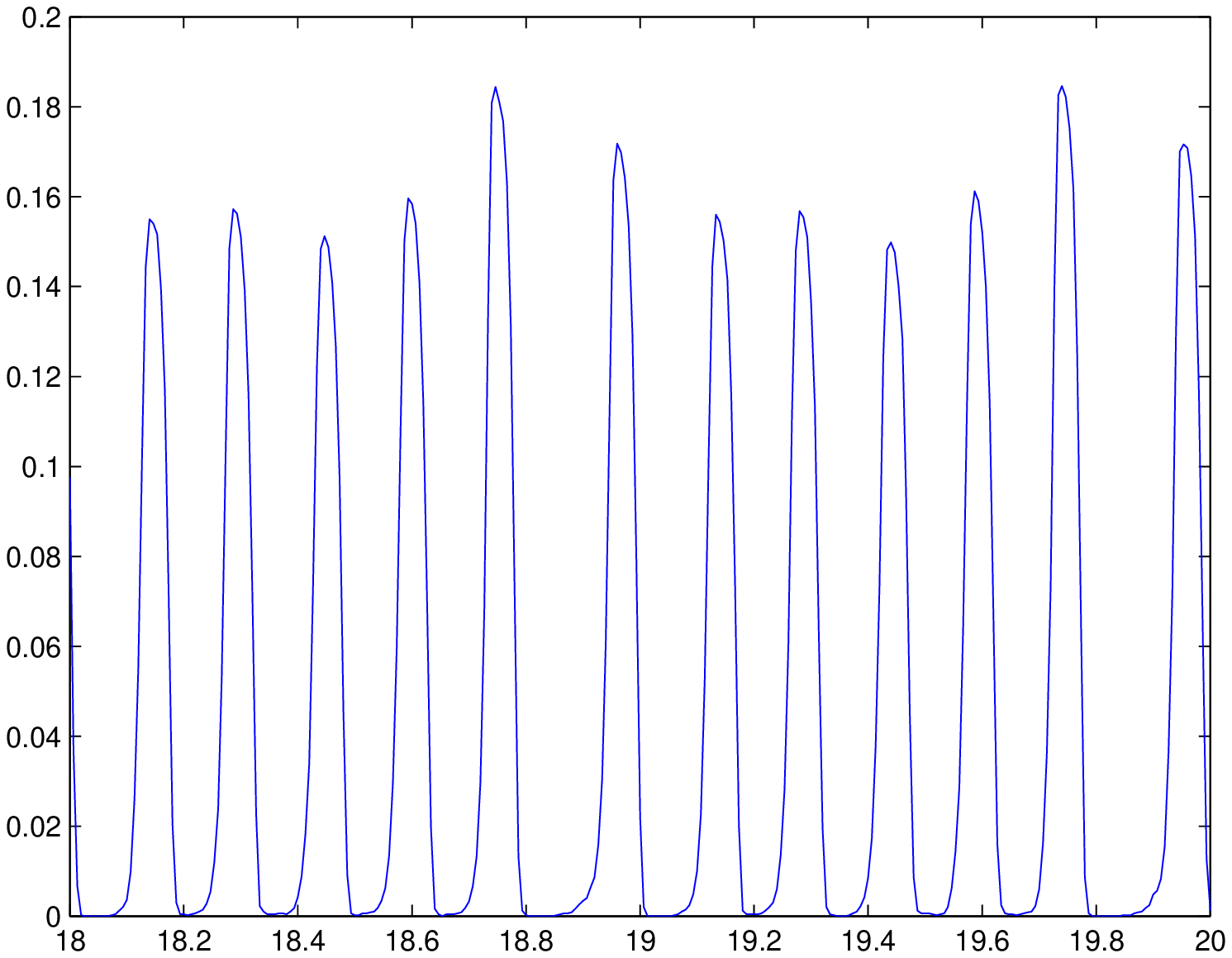,height=5cm,width=5cm}}}}
\caption{Simulations of the graduated model with positive feedback. Here
$R=[.15,.2]$, $S=[.2,.25]$, $\sigma = .01$ and the
linear factor was 10. (a) Histogram of the final distribution of
cells. (b) Time-series of the final two time frames.}
\label{figpitch3}
\end{figure}

\begin{figure}[h!]
\centerline{\hbox{{\psfig{figure=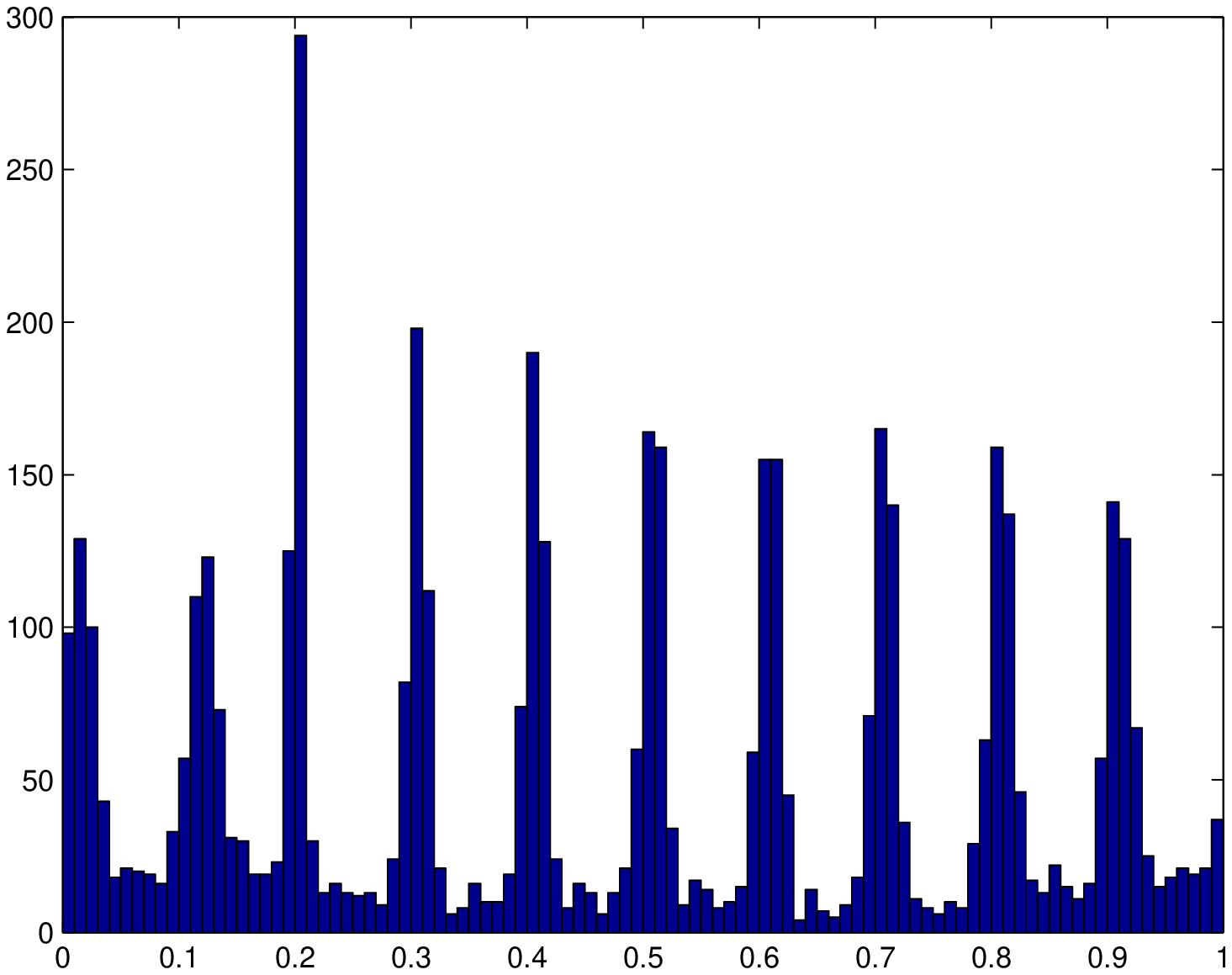,height=5cm,width=5cm}}
\hspace*{.1in}{\psfig{figure=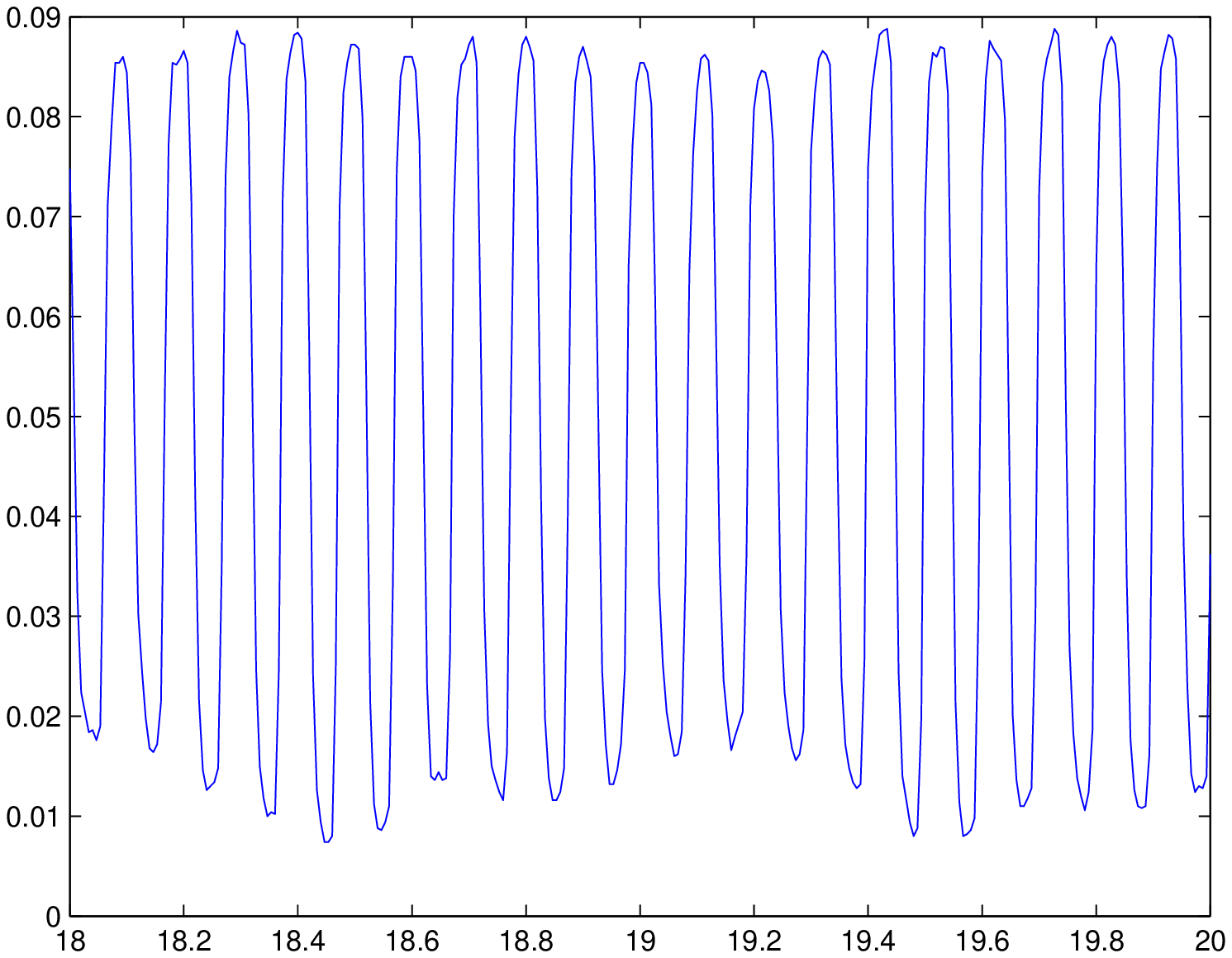,height=5cm,width=5cm}}}}
\caption{Simulations of the graduated model with negative. Here
$R=[.15,.2]$, $S=[.2,.25]$, $\sigma = .01$ and the
linear factor was $-10$. (a) Histogram of the final distribution of
cells. (b) Time-series of the final two time frames.}
\label{figpitch4}
\end{figure}

The figures clearly show clustering for all the specified parameter values.

We observe the following things from the simulations:
\begin{itemize}
 \item The negative beedback produces more clusters than positive.
 \item The number of clusters formed is inversely proportional to the width of $S$ or
       $|R|+|S|$.
 \item Positive feedback causes sharper clustering, i.e. the clusters appear to
       be completely separated.
 \item Too much noise  can destroy the clustering effect.
       The amount of noise needed to do so was smaller for larger numbers of clusters.
\end{itemize}
The first two of these observations are consistent with the results concerning
the idealized advancing and blocking models.
The third observation is somewhat in contrast to the lack of rigorous results
in section~\ref{pertadv}. The stochastically perturbed positive feedback model is seen to form
clusters just as readily and even more markedly than the negative model.
These simulations with more realistic models
reinforce the hypothesis that clustering is a robust phenomenon in cell cycle dynamics with
any form of feedback.

\subsection{Stratified Model}

When yeast cells divide, one is the mother and the other the daughter.
The mother carries a scar from the bud and the daughter does not.
The number of daughters a cell has had can be determined by counting
the bud scars. We will refer to a cell's {\em generation} as the
number of daughters it has produced, starting the count from $0$
for the daughters themselves in the first cell cycle.

It is known that after division the mother's volume is slightly
larger than that of a daughter. Further, the later generations
have slightly shorter times to budding and slightly shorter cell
cycle times.

Leslie  proposed a general stratified population model that takes
into account the differences in generations that we adapt
to yeast in \cite{tpb}.
As before, we will simplify
the model by considering logarithmic coordinates, normalized by
the cell cycle length of the zero-th generation. In these coordinates
the cell cycle of the $0$ generation is the unit interval $I_0 = [0,1]$.
We will denote the $S$ phase and $R$ phases of the $0$ generation
by $R = [r_0,s_0]$ and $S = [s_0,t_0]$. The successive generations
in these coordinates are represented by the intervals $I_k = [0,D_k]$
with $R_k = [r_k,s_k]$ and $S_k = [s_k,t_k]$.

Note first that if $D_k =1$, $R_k = R_0$ and $S_k = S_0$ for all $k\ge 1$,
then this model reduces to the simplified model of the previous sections.
Next we note that if these conditions approximately hold then the above
analysis can be applied in the same way as the small random perturbations.

If there are marked differences in the parameters between generations, then
we still might be able to repeat a rigorous analysis. For example, suppose
that $D_k = 2/3$ for all $k>0$. Then in the idealized models we could have
3 clusters in the 0-th generation and 2 clusters in the higher generations.
In this scenario, at division parents would be synchronized with a cluster
of daughters that was previously 1/3 of a cycle ahead of them.
As another example, suppose that the ratio above is 3/4, then the system
could support 4 daughter clusters and 3 parent clusters.
In the same way, rational ratios of cell cycles could produce a variety of
clustering combinations.

In real systems, first of all, ratios close to the above idealized ``resonances''
could still lead to clustering. Secondly, we may suppose that only
the first few generations matter
since higher generations are represented in numbers that decay approximately
geometrically. Biologically, the successive generations beyond the second
have not been found to exhibit marked differences anyway.

\section{Discussion}

The main conclusion we wish to emphasize is that clustering seems to be
a very robust phenomenon; it occurs in all the models studied and
for a large spectrum of parameter values. This later point is quite
important since the actual systems in question are so
complex that many of the parameter values are difficult to accurately
determine.

The observed synchronous behavior here is not driven by the cell
cycle itself, but by feedback mechanisms acting on the cell cycle.
However, the clusters must necessarily be an integer in number
and so the oscillations produced by clustering would naturally
appear with a period that is an integer fraction of the cell cycle period.

We observe two phenomena that possibly make negative feedback a more
reasonable explanation of yeast cell-cycle synchrony. First, negative feedback
allows for large numbers of clusters as in some experiments.
Second, negative feedback is seen to initiate
clustering more naturally, which we observe in analysis of the
idealized models. However, either positive or negative
feedback are possible agents of clustering and conclusive
evidence of either cause would need extensive biological modeling
confirmed by experiment. Currently there is a strong interest in
modeling the details of the cell cycle.

Given the robustness of clustering, we suspect that it
occurs in many other biological systems with cell cycles.
Specifically it could play a role in many types of microbiological
systems with cell cycles and some type of signaling, including
bacteria \cite{chen,collier,dunny,lyon,mac}

\appendix

\section{PDE models of the cell cycle}

In this appendix we turn our attention briefly
to partial differential equation approximations of the models.
The PDE models derived  are presented to give some context to the
ODE, RDE and SDE models that we consider in the rest of the paper.
Also, there is some history of such models in the literature.

\subsection{General PDE model}

Let $U(x,t)$ denote the distribution of cells within the cell cycle as represented
by $x \in [0,1)$, i.e.
$$
   U(x,t) = \frac{1}{N} \sum_{i \in S(t)} \delta_{x_i(t)}(x)
$$
where $S(t) \subset S$ is the index set of cells that are living at time $t$
and $N$ might be taken as the average number of active cells.
As is customary, we may approximate the point distribution by a more regular
distribution $u(x,t)$.

If there is no death or harvesting, then formally a distribution
$u(t,x)$ evolves locally under (\ref{sde}) by
the Fokker-Planck equation, which takes the form:
\begin{equation}\label{Fokker-Plank}
  \frac{\partial u}{\partial t} +  \frac{\partial}{\partial x} \left(a(x,t,[u])\, u\right)
          - \sigma \frac{\partial^2 u}{\partial x^2} = 0.
\end{equation}
In this notation $[u]$ denotes a functional dependence on $u$
and possibly its history, for example via an integral operator.
Note that this derivation is not standard since the
cells in the distribution are not independent from each other,
but coupled through the term $a$ and the resulting equation is
inherently nonlinear. Thus for rigour, equation (\ref{Fokker-Plank}) requires
further justification.

If the culture is well mixed and harvesting is via removal of
bulk material, as in a bio-reactor, then all cells are equally likely
to be harvested, regardless of state and the effect of harvesting
on the density will be proportional to $u(t,x)$. Similarly, if
death of a cell is equally likely at any stage of the cell cycle
then the harvesting and death can be incorporated into the Fokker-Planck
equation as a term $-k u$ on the right hand side, i.e.:
\begin{equation}\label{Fokker-Plank_harvesting}
  \frac{\partial u}{\partial t} +  \frac{\partial}{\partial x} \left( a(x,t,[u])\, u \right)
          - \sigma \frac{\partial^2 u}{\partial x^2} = -ku .
\end{equation}

To take cell division into account in the PDE model, we need to require the boundary
conditions:
\begin{equation}
u(t,0) = 2 u(t,1^-), \qquad a(t,0,[u]) u(t,0) + \sigma u_x(t,0) = 2 a(t,1^-,[u]) u(t,1^-) + 2 \sigma u_x(t,1^-).
\end{equation}

Finally, in a PDE model we may represent $a$ as an
integral operator:
$$
   a(t,x,[u]) = b(t,x) + \int_0^1 k(t,x,u(z),z) \, dz.
$$

\subsection{Periodic or Near Steady State}

Next we suppose that the culture is growing in a bioreactor
near a steady state or a periodic state so that the cell division is approximately
balanced by harvesting.
Suppose also that dependence of $a$ on all variables besides $x$ is small, or,
that the dependence on $t$ is averaged over one period,
and that $\gamma$ is small.
Then equation (\ref{Fokker-Plank_harvesting})
may be approximated by:
\begin{equation}\label{pde_harvesting}
 \frac{\partial u}{\partial t} +  \frac{\partial}{\partial x} \left( a(x)\, u \right)
           = -ku .
\end{equation}
If we make the change of variables:
$$
 u(t,x) = e^{-k\int_0^x \frac{1}{a(y)} \, dy} z(t,x),
$$
then an easy calculation shows that $z(t,x)$ satisfies:
\begin{equation}\label{z}
  \frac{\partial z}{\partial t} +  \frac{\partial}{\partial x} \left( a z \right) = 0.
\end{equation}

Further, the balance between growth and harvesting implies $k = \ln 2/(\int_0^1 \frac{1}{a})$.
This implies that the boundary conditions on $z$ are
\begin{equation}\label{clbc}
   z(t,0) = z(t,1).
\end{equation}
The ordinary differential equation that generates (\ref{z}) is simply:
$$
   \frac{dx}{dt} = a(x),
$$
on $[0,1)$.
If we reincorporate the other influences on the growth rate
then the unperturbed and pertubed equations become:
\begin{equation}\label{ode2}
   \frac{dx_i}{dt} = a(t,x_i,\bar{x}),
\end{equation}
and
\begin{equation}\label{rde2}
   \frac{dx_i}{dt} = a(t,x_i,\bar{x}) + g_i(t).
\end{equation}
Note that (\ref{z}) with (\ref{clbc}) is a conservative equation in
the sense that $\int z \, dx$ is preserved, thus in
the ODE (\ref{ode2}) and RDE (\ref{rde2})  the total number of cells considered should not
change. It is thus reasonable to take $S$ to be a fixed
set $\{1,2, \ldots, N \}$ and let $x_i(t)$ to be defined
for all $t \in [0,T]$. To accommodate division we may reset $x_i(t)$
to $0$ each time it reaches $1$.

\subsection{Notes the PDE models and related work}

The importance of the cell cycle has been recognized in
hematopoiesis~\cite{mackey07,mackey06}, as well as yeast growth~\cite{gage},
where PDE models similar to~(\ref{pde_harvesting}) are considered
with $a$ assumed constant. Thirty years ago Rotenberg~\cite{rotenberg}
described the PDE model~(\ref{pde_harvesting}) for yeast growth, but with
constant coefficients and an externally applied periodic blocking at cell division.
He demonstrated numerically that clustering occurs in the model. Surprisingly,
this article has received few citations and little attention. A trivial, but
important step, common to our work and that of Rotenberg, is normalization of
the cell cycle. This provides a standard domain on which to consider both ODE and PDE models.

In related recent work~\cite{zhu} Zhu et.~al.~have considered population
balance models coupled with a substrate balance equation, where $a = a(s)$ with
$s$ the effective substrate concentration. In~\cite{henson} $a = a(s,x)$ is
admitted, but the focus is on a survey of possible numerical methods,
not analysis of the model, and clustering solutions are not discussed.
Note that the models considered in the current manuscript implicitly
incorporate the substrate balance equation (and also the density of any
other signaling compounds) into the term $a$ via $[u]$.

\section*{Acknowledgment}

E.M.B.\ and C.C.S.\ were partially supported through NSF-DMS 0443855, NSF-ECS 0601528 and
the W.M.~Keck Foundation Grant \#062014. T.G.\ was partially supported by
NSF/NIH grant 0443863, NIH-NCRR grant PR16445  and NSF-CRCNS grant 0515290.

We would like to acknowledge the input of the late Robert Klevecz in this manuscript.
He is responsible for attracting our attention to this field, he generously
opened his lab and shared his knowledge with us. Without his support and encouragement
this manuscript would not be possible and we dedicate it to his memory.

\end{document}